\documentclass[reqno]{amsart}
\usepackage{amssymb}
\usepackage{hyperref}
\usepackage[english]{babel}

\numberwithin{equation}{section}
\newtheorem{theorem}{Theorem}[section]

\newtheorem{corollary}[theorem]{Corollary}

\newcommand{\lr}[1]{\left(#1\right)}

\begin{document}
\title[A short survey of recent results on ...]
      {A short survey of recent results on Buschman--Erd\'{e}lyi transmutations}

\author[S.M. Sitnik]
{Sergei M. Sitnik
 \\   \\  \\  \vspace*{1.3cm}
\hfill{\scriptsize
{\it Dedicated to Professor Ivan Dimovski's contributions}}\\
\vspace*{-1.5cm}}

\address{S.M. Sitnik \newline
Voronezh Institute of the Ministry of Internal Affairs, Pr. Patriotov, 53 \hfill \break
Voronezh, 394065, Russia \newline
\textit{and}\newline
 RUDN 
 University,  6 Miklukho--Maklaya Str. \hfill \break
  Moscow, 117198, Russia}
\email{pochtasms@gmail.com}

\subjclass[2000]{26A33, 44A15}
\keywords{Transmutations; Buschman--Erd\'{e}lyi transmutations; Sonine--Poisson--Delsarte transmutations,
Sonine--Katrakhov and Poisson--Katrakhov transmutations; Hardy operator; Kipriyanov space; fractional Integrals; Bessel operator}

\begin{abstract}

This short survey paper contains brief historical information, main known facts and original author's results on the theory of the Buschman--Erd\'{e}lyi transmutations and some of their applications. The operators of  Buschman--Erd\'{e}lyi  type were first studied by E.T.~Copson, R.G.~Buschman and A.~Erd\'{e}lyi  as integral operators. In 1990's the author was first to prove the transmutational nature of these operators and published papers with detailed study of their properties. This class include as special cases such famous objects as the Sonine--Poisson--Delsarte transmutations and the fractional Riemann--Lioville integrals. In this paper, the Buschman--Erdelyi transmutations are logically classified as operators of the first kind with special case of zero order smoothness operators, second kind and third kind  with special case of unitary Sonine--Katrakhov and Poisson--Katrakhov transmutations. We study such properties as transmutational conditions, factorizations, norm estimates, connections with classical integral transforms. Applications are considered to singular partial differential equations, embedding theorems with sharp constants in Kipriyanov spaces, Euler--Poisson--Darboux equation including Copson lemma, generalized translations, Dunkl operators,  Radon transform, generalized harmonics theory, Hardy operators, V.\,Katrakhov's results on pseudodifferential operators  and boundary--value problems of new kind  for equations with solutions of arbitrary growth at the isolated singularity for elliptic partial differential equations.

\end{abstract}

\maketitle

\vspace*{-0.8cm}
\section{Introduction: an idea of transmutations, \break historical information and applications}

\textit{This paper was published in \cite{Si} in the special issue of the Journal of Inequalities and Special Functions
dedicated to  Professor Ivan Dimovski's contributions to different fields of mathematics: transmutation theory, special functions, integral transforms, function theory etc.}

Let us start with the main definition.

\smallskip

\noindent
\textbf{Definition 1.}
For a given pair of operators $(A,B)$ an operator $T$ is called transmutation (or intertwining) operator if on  elements of some functional spaces the following property is valid
\begin{equation}
\label{1.1}
\Large {T\,A=B\,T.}
\end{equation}

It is obvious that the notion of transmutation is direct and far going generalization of the matrix similarity from linear algebra. But  the transmutations do not reduce to similar operators because intertwining operators often are not bounded in classical spaces and the inverse operator may  not exist or not be bounded in the same space. As a consequence, spectra of intertwining operators are not the same as a rule. Moreover transmutations may be unbounded. It is the case for the Darboux transformations which are defined for a pair of differential operators and are differential operators themselves, in this case all three operators are unbounded in classical spaces. But the theory of the Darboux transformations is included in transmutation theory too. Also a pair of intertwining operators may not be differential ones. In transmutation theory there are problems for the following various types of operators: integral, integro--differential, difference--differential (e.g.  the Dunkl operator), differential or integro--differential of infinite order (e.g. in connection with Schur's lemma), general linear operators in functional spaces, pseudodifferential and abstract  differential operators.

All classical integral transforms due to Definition 1 are also special cases of transmutations, they include the Fourier, Petzval (Laplace), Mellin, Hankel, Weierstrass, Kontorovich--Lebedev, Meyer, Stankovic, Obrechkoff, finite Grinberg and other transforms.

In quantum physics, in study of Shr\"{o}dinger equation and inverse scattering theory, the underlying transmutations are called wave operators.

The commuting operators are also a special class of transmutations. The most important class consists of operators commuting with derivatives. In this case transmutations as commutants are usually in the form of formal series, pseudodifferential or infinite order differential operators. Finding of commutants is directly connected with finding all  transmutations in the given functional space. For these problems works a  theory of operator convolutions, including the Berg--Dimovski convolutions \cite{Dim}. Also, more and more applications are developed connected with the transmutation theory for commuting differential operators, such problems are based on classical results of J.L.\,Burchnall, T.W.\,Chaundy. The transmutations are also connected with factorization problems for integral and differential operators. Special class of transmutations are the so called Dirichlet--to--Neumann and  Neumann--to--Dirichlet operators which link together solutions of the same equation but with different kinds of boundary conditions.

And how the transmutations usually works? Suppose we study properties for a rather complicated operator $A$. But suppose also that we know the corresponding  properties for a  model more simple operator $B$ and transmutation (\ref{1.1}) readily exists. Then we usually may copy results for the model operator $B$ to corresponding ones for the more complicated operator $A$. This is shortly the main idea of transmutations.

Let us  consider for example an equation $Au=f$, then applying to it a transmutation with property (\ref{1.1}) we consider a new equation $Bv=g$, with $v=Tu, g=Tf$. So if we can solve the simpler  equation $Bv=g$, then the initial one is also solved and has solution $u=T^{-1}v$. Of course, it is supposed that the inverse operator exists and its explicit form is known. This is a simple application of the transmutation technique for finding and proving formulas for solutions of ordinary and partial differential equations.

The monographs  \cite{Car1}-\cite{Tri} are completely devoted to the transmutation theory and its applications, note also author's survey \cite{Sita}. Moreover, essential parts of  monographs \cite{Dim}, \cite{CarSho}-\cite{Hro}, \cite{Kir1}, etc. include material on transmutations, the complete list of books which investigate some transmutational problems is now near of 100 items.

We use the term ``transmutation"  due to \cite{Car3}: ``Such operators are often called transformation operators by the Russian school (Levitan, Naimark, Marchenko et. al.), but transformation seems too broad term, and, since some of the machinery seems ``magical"\,  at times, we have followed Lions and Delsarte in using the word ``\textit{transmutation}''.

Now the transmutation theory is a completely formed part of the mathematical world in which methods and ideas from different areas are used: differential and integral equations, functional analysis, function theory, complex analysis, special functions, fractional integrodifferentiation.

The transmutation theory is deeply connected with many applications in different fields of mathematics. The transmutation operators are applied in inverse problems via the generalized Fourier transform, spectral function and famous Levitan equation;
in scattering theory the Marchenko equation is formulated in terms of transmutations; in spectral theory transmutations help to prove trace formulas and asymptotics for spectral function; estimates for transmutational kernels control stability in inverse and scattering problems;
for nonlinear equations via Lax method transmutations for Sturm--Lioville problems lead to proving existence and explicit formulas for soliton solutions. Special kinds of transmutations are the generalized analytic functions, generalized translations and convolutions, Darboux transformations. In the theory of partial differential equations the transmutations works for proving explicit correspondence formulas among solutions of perturbed and non--perturbed equations, for singular and degenerate equations, pseudodifferential operators, problems with essential singularities at inner or corner points, estimates of solution decay for elliptic and ultraelliptic equations. In function theory transmutations are applied to embedding theorems and generalizations of Hardy operators, Paley--Wiener theory, generalizations of harmonic analysis based on generalized translations. Methods of transmutations are used in many applied problems: investigation of Jost solutions in scattering theory, inverse problems, Dirac and other matrix systems of differential equations, integral equations with special function kernels, probability theory and random processes, stochastic random equations, linear  stochastic estimation, inverse problems of geophysics and transsound gas dynamics. Also a number of applications of the transmutations to nonlinear equations is permanently increased.

 In fact, the modern transmutation theory originated from two basic examples, see \cite{Sita}.
The first is the transmutation $T$ for Sturm--Lioville problems with some potential $q(x)$ and natural boundary conditions
\vskip -10pt
\begin{equation}
T(D^2\,y(x)+q(x)y(x))=D^2\,(Ty(x)), D^2\,y(x)=y''(x),
\end{equation}
The second example is a problem of studying transmutations intertwining the Bessel operator $B_\nu$ and the second derivative:
\vskip -10pt
\begin{equation}
T\lr{B_\nu} f=\lr{D^2} Tf,\ \  B_{\nu}=D^2+\frac{2\nu +1}{x}D,\ \  D^2=\frac{d^2}{dx^2},\ \  \nu \in \mathbb{C}.
\end{equation}
This class of transmutations includes the Sonine--Poisson--Delsarte, Buschman--Erd\'{e}lyi  operators and their generalizations. Such transmutations found many applications for a special class of partial differential equations with singular coefficients. A typical equation of this class is the $B$--elliptic equation with the Bessel operator in some variables of the form
\vskip -10pt
\begin{equation}
\label{59}
\sum_{k=1}^{n}B_{\nu,x_k}u(x_1,\dots, x_n)=f.
\end{equation}
Analogously,  $B$--hyperbolic and $B$--parabolic equations are considered, this terminology was proposed by I.\,Kipriyanov. This class of equations was first studied by Euler, Poisson, Darboux and continued in Weinstein's theory of generalized axially symmetric potential (GASPT). These problems were further investigated by Zhitomirslii, Kudryavtsev, Lizorkin, Matiychuk, Mikhailov, Olevskii, Smirnov, Tersenov, He Kan Cher, Yanushauskas, Egorov and others.

In the most detailed and complete way, equations with Bessel operators were studied by the Voronezh mathematician Kipriyanov and his disciples Ivanov, Ryzhkov, Katrakhov, Arhipov, Baidakov, Bogachov, Brodskii, Vinogradova, Zaitsev, Zasorin, Kagan, Katrakhova, Kipriyanova, Kononenko, Kluchantsev, Kulikov, Larin, Leizin, Lyakhov, Muravnik, Polovinkin, Sazonov, Sitnik, Shatskii, Yaroslavtseva. The essence of Kipriyanov's school results was published in  \cite{Kip1}. For classes of equations with Bessel operators, Kipriyanov introduced special functional spaces which were named after him \cite{Kip2}.
In this field interesting results were investigated by Katrakhov and his disciples, now these problems are considered by  Gadjiev, Guliev, Glushak, Lyakhov, Shishkina with their coauthors and students. Abstract equations of the form (\ref{59}) originated from the monograph \cite{CarSho} were considered by  Egorov, Repnikov, Kononenko, Glushak, Shmulevich and others. And transmutations are one of basic tools for equations with Bessel operators, they are applied to construction of solutions, fundamental solutions, study of singularities, new boundary--value and other problems.

We must note that the term ``operator"\, is used in this paper for brevity in the broad and sometimes not exact meaning, so appropriate domains and function classes are not always specified. It is easy to complete and make strict for every special result.

\vspace*{-0.45cm}

\section{Buschman--Erd\'{e}lyi  transmutations}

The term ``Buschman--Erd\'{e}lyi  transmutations"\, was introduced by the author and is now accepted. Integral
equations with these operators were studied in mid--1950th. The author was first to prove the transmutational nature of these operators. The classical Sonine and Poisson operators are special cases of the Buschman--Erd\'{e}lyi  transmutations and Sonine--Dimovski and Poisson--Dimovski transmutations are their generalizations for the hyper--Bessel equations and functions.

The Buschman--Erd\'{e}lyi transmutations have many modifications. The author introduced convenient classification of them. Due to this classification we introduce Buschman--Erd\'{e}lyi  transmutations of the first kind, their kernels are expressed in terms of Legendre functions of the first kind. In the limiting case we define Buschman--Erd\'{e}lyi  transmutations of zero order smoothness being important in applications. The kernels of Buschman--Erd\'{e}lyi  transmutations of the second  kind are expressed in terms of Legendre functions of the second kind. Some combination of operators of  the first kind and the second kind leads to  operators of the third kind. For the special choice of parameters they are unitary operators in the standard Lebesgue space.  The author proposed the terms
``Sonine--Katrakhov"\, and ``Poisson--Katrakhov"\, transmutations in honor of V.\,Katrakhov who  introduced and studied  these operators.

The study of integral equations and invertibility for the Buschman--Erd\'{e}lyi  operators  was started in 1960-th by P. Buschman and A. Erd\'{e}lyi,  \cite{Bu1}--\cite{Er2}. These operators also were investigated by Higgins, Ta Li, Love, Habibullah, K.N.\,Srivastava, Ding Hoang An, Smirnov, Virchenko, Fedotova, Kilbas, Skoromnik and others.
During this period, for this class of operators were considered only problems of solving integral equations, factorization and invertibility, cf.  \cite{KK}.

The most detailed study of the Buschman--Erd\'{e}lyi  transmutations was taken by the author in 1980--1990th
\cite{Sit3}--\cite{Sit4} and continued in \cite{Sit2}--\cite{Sit14} and some other papers.
Interesting and important results were proved by N. Virchenko and A. Kilbas and their disciples \cite{KiSk1}--\cite{ViFe}.

Let us first consider  the most well--known transmutations for the Bessel operator and the second derivative:
\vskip -12pt
\begin{equation}
\label{2.1}
T\lr{B_\nu} f=\lr{D^2} Tf, \ B_{\nu}=D^2+\frac{2\nu +1}{x}D, \ D^2=\frac{d^2}{dx^2}, \ \, \nu \in \mathbb{C}.
\end{equation}


\noindent
\textbf{Definition 2.}
The Poisson transmutation is defined by
\vskip -10pt
\begin{equation}
\label{2.2}
P_{\nu}f=\frac{1}{\Gamma(\nu+1)2^{\nu}x^{2\nu}}
\int_0^x \left( x^2-t^2\right)^{\nu-\frac{1}{2}}f(t)\,dt,\ \ \Re \nu> -\frac{1}{2}.
\end{equation}
Respectivelt, the Sonine transmutation is defined by
\vskip -12pt
\begin{equation}
\label{2.3}
S_{\nu}f=\frac{2^{\nu+\frac{1}{2}}}{\Gamma(\frac{1}{2}-\nu)}\frac{d}{dx}
\int_0^x \left( x^2-t^2\right)^{-\nu-\frac{1}{2}}t^{2\nu+1}f(t)\,dt,\ \  \Re \nu< \frac{1}{2}.
\end{equation}
The operators (\ref{2.2})--(\ref{2.3}) intertwine by the formulas
\vskip -10pt
\begin{equation}
\label{56}
S_\nu B_\nu=D^2 S_\nu,\ \  P_\nu D^2=B_\nu P_\nu.
\end{equation}
The definition may be extended to $\nu\in\mathbb{C}$.
We will use more historically exact term as the Sonine--Poisson--Delsarte transmutations, \cite{Sita}.

An important generalization for the Sonine--Poisson--Delsarte are the transmutations for the hyper--Bessel operators and functions.
Such functions were first considered by Kummer and Delerue. The detailed study on these operators and hyper--Bessel functions was done
by Dimovski and further, by  Kiryakova. 
The corresponding  transmutations have been called by Kiryakova \cite{Kir1} as the Sonine--Dimovski and  Poisson--Dimovski transmutations.
In hyper--Bessel operators theory the leading role is for the Obrechkoff integral transform \cite{Dim1}--\cite{DimKir1}, \cite{Kir1}. It is a transform with Meijer's $G$--function kernel which generalizes the Laplace, Meijer and many other integral transforms introduced by different authors. Various results on the hyper--Bessel functions, connected equations and transmutations were many times reopened. The same is true for the Obrechkoff integral transform. It my opinion, the Obrechkoff  transform together  with the Laplace, Fourier, Mellin, Stankovic transforms are essential basic elements from which many other transforms are constructed with corresponding applications.

Let us define and study some main properties of  the Buschman--Erd\'{e}lyi  transmutations of the first kind.
This class of transmutations for some choice of parameters generalize the Sonine--Poisson--Delsart  transmutations, Riemann--Liouville  and Erd\'{e}lyi--Kober fractional integrals, Mehler--Fock transform.

\smallskip

\noindent
\textbf{Definition 3.}
Define the Buschman--Erd\'{e}lyi  operators of the first kind by
\begin{eqnarray}
\label{71}
B_{0+}^{\nu,\mu}f=\int_0^x \left( x^2-t^2\right)^{-\frac{\mu}{2}}P_\nu^\mu \left(\frac{x}{t}\right)f(t)d\,t,\\
E_{0+}^{\nu,\mu}f=\int_0^x \left( x^2-t^2\right)^{-\frac{\mu}{2}}\mathbb{P}_\nu^\mu \left(\frac{t}{x}\right)f(t)d\,t,\\
B_{-}^{\nu,\mu}f=\int_x^\infty \left( t^2-x^2\right)^{-\frac{\mu}{2}}P_\nu^\mu \left(\frac{t}{x}\right)f(t)d\,t,\\
\label{72}
E_{-}^{\nu,\mu}f=\int_x^\infty \left( t^2-x^2\right)^{-\frac{\mu}{2}}\mathbb{P}_\nu^\mu \left(\frac{x}{t}\right)f(t)d\,t.\\ \nonumber
\end{eqnarray}
Here $P_\nu^\mu(z)$ is the Legendre function of the first kind, $\mathbb{P}_\nu^\mu(z)$ is this function on the cut $-1\leq  t \leq 1$ (\cite{BE}), $f(x)$  is a locally summable function with some growth conditions at $x\to 0,\, x\to\infty$. The parameters are $\mu,\nu\in\mathbb{C}$, $\Re \mu <1$,  $\Re \nu \geq -1/2$.

Now consider some main properties for this class of transmutations, following essentially \cite{Sit3}, \cite{Sit4}, and also \cite{Sita}, \cite{Sit2}. All functions further are defined on positive semiaxis.
So we use notations $L_2$ for the functional space $L_2(0, \infty)$ and $L_{2, k}$ for power weighted
space $L_{2, k}(0, \infty)$ equipped with norm
\begin{equation}
\int_0^\infty |f(x)|^2 x^{2k+1}\,dx,
\end{equation}
$\mathbb{N}$ denotes the set of naturals, $\mathbb{N}_0$--positive integer, $\mathbb{Z}$--integer and  $\mathbb{R}$--real numbers.

First, add to Definition 3 a case of parameter $\mu =1$. It defines a very important class of operators.

\smallskip

\noindent
\textbf{Definition 4.}
Define for $\mu =1$ the Buschman--Erd\'{e}lyi  operators of zero order smoothness by
\vskip -10pt
\begin{eqnarray}
\label{2BE01}
B_{0+}^{\nu,1}f={_1 S^{\nu}_{0+}f}=\frac{d}{dx}\int_0^x P_\nu \left(\frac{x}{t}\right)f(t)\,dt,\\
\label{2BE02}
E_{0+}^{\nu,1}f={_1 P^{\nu}_{-}}f=\int_0^x P_\nu \left(\frac{t}{x}\right)\frac{df(t)}{dt}\,dt,\\
\label{2BE03}
B_{-}^{\nu,1}f={_1 S^{\nu}_{-}}f=\int_x^\infty P_\nu \left(\frac{t}{x}\right)(-\frac{df(t)}{dt})\,dt,\\
\label{2BE04}
E_{-}^{\nu,1}f={_1 P^{\nu}_{0+}}f=(-\frac{d}{dx})\int_x^\infty P_\nu \left(\frac{x}{t}\right)f(t)\,dt,
\end{eqnarray}
where $P_\nu(z)=P_\nu^0(z)$ is the Legendre function.

\begin{theorem}
The next formulas hold true for factorizations of Buschman--Erd\'{e}lyi  transmutations for suitable functions via Riemann--Liouville fractional integrals and  Buschman--Erd\'{e}lyi  operators of zero order smoothness:
 \begin{equation}\label{1.9}
{B_{0+}^{\nu,\,\mu} f=I_{0+}^{1-\mu}~ {_1 S^{\nu}_{0+}f},~B_{-}^{\nu, \,\mu} f={_1 P^{\nu}_{-}}~ I_{-}^{1-\mu}f,}
\end{equation}
\begin{equation}\label{1.10}
{E_{0+}^{\nu,\,\mu} f={_1 P^{\nu}_{0+}}~I_{0+}^{1-\mu}f,~E_{-}^{\nu, \, \mu} f= I_{-}^{1-\mu}~{_1 S^{\nu}_{-}}f.}
\end{equation}
\end{theorem} 

These formulas allow to separate parameters $\nu$ and $\mu$. We will prove soon that operators
\eqref{2BE01}--\eqref{2BE04} are isomorphisms of $L_2(0, \infty)$ except for some special parameters. So,
 operators  \eqref{71}--\eqref{72} roughly speaking are of the same smoothness in $L_2$ as integrodifferentiations $I^{1-\mu}$ and they coincide with them for $\nu=0$. It is also possible to define Buschman--Erd\'{e}lyi  operators for all $\mu\in\mathbb{C}$.

\smallskip

\noindent
\textbf{Definition 5.}
Define the number $\rho=1-Re\,\mu $ as smoothness order for Buschman--Erd\'{e}lyi  operators \eqref{71}--\eqref{72}.

\smallskip

So for $\rho > 0$ (otherwise for $Re\, \mu > 1$) the Buschman--Erd\'{e}lyi  operators are smoothing and for $\rho < 0$ (otherwise for $Re\, \mu < 1$) they decrease smoothness in $L_2$ spaces. Operators \eqref{2BE01}--\eqref{2BE04} for which $\rho = 0$ due to Definition 5 are of zero smoothness order in accordance with their definition.

For some special parameters $\nu,~\mu$ the Buschman--Erd\'{e}lyi  operators of the first kind are reduced to other known operators. So for $\mu=-\nu$ or $\mu=\nu+2$ they reduce to Erd\'{e}lyi--Kober operators, for $\nu = 0$ they reduce to fractional integrodifferentiation $I_{0+}^{1-\mu}$ or $I_{-}^{1-\mu}$, for  $\nu=-\frac{1}{2}$, $\mu=0$ or $\mu=1$ kernels reduce to elliptic integrals, for $\mu=0$,  $x=1$, $v=it-\frac{1}{2}$ the operator $B_{-}^{\nu, \, 0}$ differs only by a constant from Mehler--Fock transform.

As a pair for the Bessel operator consider a connected one
\begin{equation}
\label{75}
L_{\nu}=D^2-\frac{\nu(\nu+1)}{x^2}=\left(\frac{d}{dx}-\frac{\nu}{x}\right)
\left(\frac{d}{dx}+\frac{\nu}{x}\right),
\end{equation}
which for $\nu \in \mathbb{N}$ is an angular momentum operator from quantum physics. Their transmutational relations are established in the next theorem.

\begin{theorem}
For a given pair of transmutations $X_\nu, Y_\nu$
\begin{equation}
\label{76}
X_\nu L_{\nu}=D^2 X_\nu , Y_\nu D^2 = L_{\nu} Y_\nu
\end{equation}
define the new pair of transmutations by formulas
\begin{equation}
\label{77}
S_\nu=X_{\nu-1/2} x^{\nu+1/2}, P_\nu=x^{-(\nu+1/2)} Y_{\nu-1/2}.
\end{equation}
Then for the new pair  $S_\nu, P_\nu$ the next formulas are valid:
\begin{equation}
\label{78}
S_\nu B_\nu = D^2 S_\nu, P_\nu D^2 = B_\nu P_\nu.
\end{equation}
\end{theorem}

\begin{theorem}\label{t6}
Let $Re \, \mu \leq 1$. Then an operator $B_{0+}^{\nu, \, \mu}$ on proper functions is a  Sonine type transmutation and
\eqref{76} is valid.
\end{theorem}

The same result holds true for other Buschman--Erd\'{e}lyi  operators, $E_{-}^{\nu, \, \mu}$ is Sonine type
and $E_{0+}^{\nu, \, \mu}$,  $B_{-}^{\nu, \, \mu}$ are Poisson type transmutations.

From these transmutation connections, we conclude that  the Buschman--Erd\'{e}lyi  operators link the corresponding eigenfunctions for the two operators. They lead to formulas for the Bessel functions via exponents and trigonometric functions, and vice versa which generalize the classical Sonine and Poisson formulas.

Now consider factorizations of the Buschman--Erd\'{e}lyi  operators.
First let us list the main forms of fractional integrodifferentiations:
Riemann--Liouville, Erd\'{e}lyi--Kober, fractional integral by function $g(x)$, cf. \cite{KK},
\begin{eqnarray}
\label{61}
I_{0+,x}^{\alpha}f=\frac{1}{\Gamma(\alpha)}\int_0^x \left( x-t\right)^{\alpha-1}f(t)d\,t,\\ \nonumber
I_{-,x}^{\alpha}f=\frac{1}{\Gamma(\alpha)}\int_x^\infty \left( t-x\right)^{\alpha-1}f(t)d\,t,
\end{eqnarray}
\begin{eqnarray}
\label{62}
I_{0+,2,\eta}^{\alpha}f=\frac{2 x^{-2\lr{\alpha+\eta}}}{\Gamma(\alpha)}\int_0^x \left( x^2-t^2\right)^{\alpha-1}t^{2\eta+1}f(t)d\,t,\\ \nonumber
I_{-,2,\eta}^{\alpha}f=\frac{2 x^{2\eta}}{\Gamma(\alpha)}\int_x^\infty \left( t^2-x^2\right)^{\alpha-1}t^{1-2\lr{\alpha+\eta}}f(t)d\,t,
\end{eqnarray}
\begin{eqnarray}
\label{63}
I_{0+,g}^{\alpha}f=\frac{1}{\Gamma(\alpha)}\int_0^x \left( g(x)-g(t)\right)^{\alpha-1}g'(t)f(t)d\,t,\\ \nonumber
I_{-,g}^{\alpha}f=\frac{1}{\Gamma(\alpha)}\int_x^\infty \left( g(t)-g(x)\right)^{\alpha-1}g'(t)f(t)d\,t.
\end{eqnarray}
In all cases $\Re\alpha>0$ and the operators may be further defined for all $\alpha$, see \cite{KK}.
In the case of $g(x)=x$  (\ref{63}) reduces to the Riemann--Liouville integral, in the case of $g(x)=x^2$ (\ref{63}) reduces to the Erd\'{e}lyi--Kober operator, and in the case of $g(x)=\ln x$ -- to the Hadamard fractional integrals.

\begin{theorem}\label{tfact7}
The following factorization formulas are valid for the Buschman--Erd\'{e}lyi  operators
of the first kind via the Riemann--Liouville and Erd\'{e}lyi--Kober fractional integrals:
\vskip -12pt
\begin{eqnarray}
& & B_{0+}^{\nu, \, \mu}=I_{0+}^{\nu+1-\mu} I_{0+; \, 2, \, \nu+ \frac{1}{2}}^{-(\nu+1)} {\lr{\frac{2}{x}}}^{\nu+1}\label{2.17}, \\
& & E_{0+}^{\nu, \, \mu}= {\lr{\frac{x}{2}}}^{\nu+1} I_{0+; \, 2, \, - \frac{1}{2}}^{\nu+1} I_{0+}^{-(\nu+\mu)}  \label{2.18}, \\
& & B_{-}^{\nu, \, \mu}= {\lr{\frac{2}{x}}}^{\nu+1}I_{-; \, 2, \, \nu+ 1}^{-(\nu+1)} I_{-}^{\nu - \mu+2}  \label{2.19}, \\
& & E_{-}^{\nu, \, \mu}= I_{-}^{-(\nu+\mu)} I_{-; \, 2, \, 0} ^{\nu+1} {\lr{\frac{x}{2}}}^{\nu+1}  \label{2.20}.
\end{eqnarray}
\end{theorem}

The Sonine--Poisson--Delsarte transmutations also are special cases for this class of operators.

Now let us study the properties of the Buschman--Erd\'{e}lyi  operators of zero order smoothness, defined by \eqref{2BE01}--\eqref{2BE04}.
A similar operator was introduced by Katrakhov by multiplying the Sonine operator with a fractional integral, his aim was to work with transmutation obeying good estimates in $L_2(0,\infty)$.

\vskip 2pt

We use the Mellin transform defined by \cite{Mar}
\begin{equation}
\label{710}
g(s)=M{f}(s)=\int_0^\infty x^{s-1} f(x)\,dx.
\end{equation}
\vskip -3pt \noindent
The Mellin convolution is defined by
\vskip -10pt
\begin{equation}
\label{711}
(f_1*f_2)(x)=\int_0^\infty  f_1\left(\frac{x}{y}\right) f_2(y)\,\frac{dy}{y},
\end{equation}
so the convolution operator with kernel $K$ acts under the Mellin transform as a multiplication on multiplicator
\begin{eqnarray}
\label{con}
M[Af](s)=M\ [\int_0^\infty  K\left(\frac{x}{y}\right) f(y)\,\frac{dy}{y}](s)=M[K*f](s)
=m_A(s)M{f}(s),\\ \nonumber
m_A(s)=M[K](s).\phantom{1111111111111111111}
\end{eqnarray}

We observe that the Mellin transform is a generalized Fourier transform on semiaxis with Haar measure
$\frac{dy}{y}$, \cite{Hel}. It plays important role for the theory of special functions, for example the gamma function is a Mellin transform of the exponential. With the Mellin transform the important breakthrough in evaluating integrals was done in 1970th when mainly by O. Marichev, the famous Slater's theorem was adapted for calculations. The Slater's theorem taking the Mellin transform as input gives the function itself as output via hypergeometric functions, see \cite{Mar}. This theorem occurred to be the milestone of powerful computer method for calculating integrals for many problems in differential and integral equations. The  package {\sl Mathematica} of  Wolfram Research is based on this theorem in calculating integrals.

\begin{theorem}\label{tnorm}
The Buschman--Erd\'{e}lyi  operator of zero order smoothness $_1 S^{\nu}_{0+}$ defined by (\ref{2BE01})
acts under the Mellin transform as convolution (\ref{con}) with multiplicator
\begin{equation}\label{s1}
m(s)=\frac{\Gamma(-s/2+\frac{\nu}{2}+1)\Gamma(-s/2-\frac{\nu}{2}+1/2)}
{\Gamma(1/2-\frac{s}{2})\Gamma(1-\frac{s}{2})}
\end{equation}
 for \, $\Re s<\min(2+\Re\nu,1-\Re\nu)$. Its norm is a periodic in $\nu$ and equals
\begin{equation}\label{s2}
\|B_{0+}^{\nu,1}\|_{L_2}=\frac{1}{\min(1,\sqrt{1-\sin\pi\nu})}.
\end{equation}
This operator is bounded in $L_2(0,\infty)$ if $\nu\neq 2k+1/2, k\in \mathbb{Z}$ and unbounded if
$\nu= 2k+1/2, k\in \mathbb{Z}$.
\end{theorem}

\begin{corollary}
The norms of operators \eqref{2BE01}--\eqref{2BE04} are periodic in $\nu$ with period 2  $\|X^{\nu}\|=\|X^{\nu+2}\|$,  $X^{\nu}$ is any of operators \eqref{2BE01}--\eqref{2BE04}.
\end{corollary}

\begin{corollary}
The norms of the operators ${_1 S_{0+}^{\nu}}$, ${_1 P_{-}^{\nu}}$ are not bounded in general, every norm is greater or equals to 1.
The norms are equal to 1 if $\sin \pi \nu \leq 0$.  The operators ${_1 S_{0+}^{\nu}}$, ${_1 P_{-}^{\nu}}$ are unbounded in $L_2$ if and only if $\sin \pi \nu = 1$ (or $\nu=(2k) + 1/2,~k \in \mathbb{Z}$).
\end{corollary}

\begin{corollary}
The norms of the operators ${_1 P_{0+}^{\nu}}$, ${_1 S_{-}^{\nu}}$ are all bounded in $\nu$, every norm is not greater then $\sqrt{2}$.
The norms are equal to 1 if $\sin \pi \nu \geq 0$.
The operators ${_1 P_{0+}^{\nu}}$, ${_1 S_{-}^{\nu}}$ are bounded in $L_2$ for all $\nu$.
The maximum of norm equals $\sqrt 2 $ is achieved if and only if $\sin \pi \nu = -1$ (èëè $\nu= -1/2+(2k) ,~k \in \mathbb{Z}$).
\end{corollary}

The most important property of the Buschman--Erd\'{e}lyi  operators of zero order smoothness is the unitarity for integer $\nu$. It is just the case if we interpret for these parameters the operator  $L_{\nu}$ as angular momentum operator in quantum mechanics.

\begin{theorem}
 The  operators \eqref{2BE01}--\eqref{2BE04} are unitary in  $L_2$ if and only if the parameter $\nu$ is an integer. In this case the pairs of operators
$({_1 S_{0+}^{\nu}}$, ${_1 P_{-}^{\nu}})$ and  $({_1 S_{-}^{\nu}}$, ${_1 P_{0+}^{\nu}})$
are mutually inverse.
\end{theorem}

To formulate an interesting special case, let us suppose that operators \eqref{2BE01}--\eqref{2BE04} act on functions permitting outer or inner differentiation in integrals, it is enough to suppose that $x f(x) \to 0$ for $x \to 0$. Then for  $\nu=1$
\begin{equation}\label {3.25}
{_1{P_{0+}^{1}}f=(I-H_1)f,~_1{S_{-}^{1}}f=(I-H_2)f,}
\end{equation}
and $H_1,~ H_2$ are the famous Hardy operators,
\begin{equation}\label {3.26}
{H_1 f = \frac{1}{x} \int\limits_0^x f(y) dy,~H_2 f = \int\limits_x^{\infty}  \frac{f(y)}{y} dy,}
\end{equation}
$I$ is the identic operator.

\begin{corollary}
The operators \eqref{3.25} are unitary in  $L_2$ and mutually inverse.
They are transmutations for the pair of differential operators  $d^2 / d x^2$ and $d^2 / d x^2 - 2/x^2$.
\end{corollary}

The unitarity of the shifted Hardy operators \eqref{3.25} in $L_2$ is a known fact \cite{KuPe}.
Below in application section, we introduce a new class of generalizations for the classical Hardy operators.

Now we list some properties of the operators acting as convolutions by the formula \eqref{con} and with some multiplicator under the Mellin transform and being transmutations for the second derivative and angular momentum operator in quantum mechanics.

\begin{theorem}\label{tOPmult}
Let an operator $S_{\nu}$ act by formulas \eqref{con}
and \eqref{76}. Then:

a) its multiplicator satisfies a functional equation
\vskip -10pt
\begin{equation}\label{5.1}
{m(s)=m(s-2)\frac{(s-1)(s-2)}{(s-1)(s-2)-\nu(\nu+1)};}
\end{equation}
\vskip -3pt

b) if any function $p(s)$ is periodic with period 2 $(p(s)=p(s-2))$,
 then a function $p(s)m(s)$ is a multiplicator for a new transmutation operator $S_2^{\nu}$ also acting by the rule \eqref{76}.
\end{theorem}

This theorem confirms the importance of studying transmutations in terms of the Mellin transform and
multiplicator functions.

Define the Stieltjes transform by (cf. \cite{KK})
\vskip -10pt
$$
(S f)(x)= \int\limits_0^{\infty} \frac{f(t)}{x+t} dt.
$$
This operator also acts by the formula \eqref{con} with multiplicator $p(s)= \pi /\sin (\pi s)$, it is bounded in $L_2$.
Obviously  $p(s)=p(s-2)$. So from Theorem \ref{tOPmult} it follows 
 a convolution of the Stieltjes transform with bounded transmutations \eqref{2BE01}--\eqref{2BE04}, 
  also transmutations of the same class bounded in $L_2$.

In this way many new classes of transmutations were introduced with special functions as kernels.

\vspace*{-0.3cm} 
\section{Sonine--Katrakhov and Poisson--Katrakhov transmutations}

Now we construct transmutations which are unitary for all $\nu$. They are defined by formulas
\vskip -12pt
\begin{eqnarray}
& & S_U^{\nu} f = - \sin \frac{\pi \nu}{2}\  {_2S^{\nu}}f+ \cos \frac{\pi \nu}{2}\  {_1S_-^{\nu}}f, \label{6.14} \\
& & P_U^{\nu} f = - \sin \frac{\pi \nu}{2}\  {_2P^{\nu}}f+ \cos \frac{\pi \nu}{2}\  {_1P_-^{\nu}}f. \label{6.15}
\end{eqnarray}
\vskip -3pt \noindent
For all values $\nu \in \mathbb{R}$ they are linear combinations of Buschman--Erd\'{e}lyi  transmutations of the first and second kinds of zero order smoothness. Also they are in the defined below class of Buschman--Erd\'{e}lyi  transmutations of the third kind.
The following integral representations are valid:
\begin{eqnarray}
& & S_U^{\nu} f = \cos \frac{\pi \nu}{2} \left(- \frac{d}{dx} \right) \int\limits_x^{\infty} P_{\nu}\lr{\frac{x}{y}} f(y)\,dy   \label{6.16}\\
& & + \frac{2}{\pi} \sin \frac{\pi \nu}{2} \left(  \int\limits_0^x (x^2\!-\!y^2)^{-\frac{1}{2}}Q_{\nu}^1 \lr{\frac{x}{y}} f(y)\,dy  \right.
\!-\!
 \int\limits_x^{\infty} (y^2\!-\!x^2)^{-\frac{1}{2}}\mathbb{Q}_{\nu}^1 \lr{\frac{x}{y}} f(y)\,dy \Biggl. \Biggr), \nonumber 
\end{eqnarray}
 \begin{eqnarray}
  & & P_U^{\nu} f = \cos \frac{\pi \nu}{2}  \int\limits_0^{x} P_{\nu}\lr{\frac{y}{x}} \left( \frac{d}{dy} \right) f(y)\,dy  \label{6.17} \\
& &  -\frac{2}{\pi} \sin \frac{\pi \nu}{2} \left( - \int\limits_0^x (x^2\!-\!y^2)^{-\frac{1}{2}}\mathbb{Q}_{\nu}^1\lr{\frac{y}{x}} f(y)\,dy \right.
\!-\!
\int\limits_x^{\infty} (y^2\!-\!x^2)^{-\frac{1}{2}} Q_{\nu}^1 \lr{\frac{y}{x}} f(y)\,dy \Biggl. \Biggr). \nonumber
\end{eqnarray}

\begin{theorem}
The operators \eqref{6.14}--\eqref{6.15},  \eqref{6.16}--\eqref{6.17} for all $\nu\in\mathbb{R}$
are unitary, mutually inverse and conjugate in $L_2$. They are transmutations acting by \eqref{75}.
$S_U^{\nu}$ is a Sonine type transmutation and $P_U^{\nu}$ is a Poisson type one.
\end{theorem}

Transmutations like (\ref{6.16})--(\ref{6.17}) but with kernels in more complicated form with hypergeometric functions were first introduced by Katrakhov in 1980. Due to this, the author proposed terms for this class of operators as Sonine--Katrakhov and Poisson--Katrakhov.
In author's papers these operators were reduced to more simple form of Buschman--Erd\'{e}lyi  ones.
It made possible to include this class of operators in general composition (or factorization) method  \cite{SiKa2}, \cite{SiKa3}, \cite{Sit5}.

\section{Applications of Buschman--Erd\'{e}lyi transmutations}

\subsection{Copson lemma}
\medskip

Consider the partial differential equation with two variables on the plane
$$
\frac{\partial^2 u(x,y)}{\partial x^2}+\frac{2\alpha}{x}\frac{\partial u(x,y)}{\partial x}=
\frac{\partial^2 u(x,y)}{\partial y^2}+\frac{2\beta}{y}\frac{\partial u(x,y)}{\partial y}
$$
(this is EPD equation or $B$--hyperbolic one in Kipriyanov's terminology)
for $x>0,\  y>0$ and $\beta>\alpha>0$
with boundary conditions on the characteristics
$$
u(x,0)=f(x), u(0,y)=g(y), f(0)=g(0).
$$

It is supposed that the solution $u(x,y)$ is continuously differentiable in the closed first quadrant and has second derivatives in this open  quadrant, boundary functions $f(x), g(y)$ are differentiable.

Then if the solution exists, the next formulas hold true:
\begin{equation}
\label{Cop1}
\frac{\partial u}{\partial y}=0, \ y=0, \  \frac{\partial u}{\partial x}=0, \ x=0,
\end{equation}
\begin{equation}
\label{Cop2}
2^\beta \Gamma(\beta+\frac{1}{2})\int_0^1 f(xt)t^{\alpha+\beta+1}
\lr{1-t^2}^{\frac{\beta -1}{2}}P_{-\alpha}^{1-\beta}{t}\,dt
\end{equation}
\begin{equation*}
=2^\alpha \Gamma(\alpha+\frac{1}{2})\int_0^1 g(xt)t^{\alpha+\beta+1}
\lr{1-t^2}^{\frac{\alpha -1}{2}}P_{-\beta}^{1-\alpha}{t}\,dt,
\end{equation*}
$$
\Downarrow
$$
\begin{equation}
\label{Cop3}
g(y)=\frac{2\Gamma(\beta+\frac{1}{2})}{\Gamma(\alpha+\frac{1}{2})
\Gamma(\beta-\alpha)}y^{1-2\beta}
\int_0^y x^{2\alpha-1}f(x)
\lr{y^2-x^2}^{\beta-\alpha-1}x \,dx,
\end{equation}
where  $P_\nu^\mu(z)$ is the Legendre function of the first kind \cite{Sita}.

So the main conclusion from the Copson lemma is that the data on characteristics can not be taken arbitrary,
these functions must be connected by the Buschman--Erd\'{e}lyi  operators of the first kind, for more detailed consideration cf. \cite{Sita}.

\subsection{Norm estimates and embedding theorems in Kipriyanov spaces.} 

Consider a set of functions $\mathbb{D}(0, \infty)$. If $f(x) \in \mathbb{D}(0, \infty)$ then $f(x) \in C^{\infty}(0, \infty),~ f(x)$ is zero at infinity. On this set, define the seminorms
\begin{eqnarray}
& & \|f\|_{h_2^{\alpha}}=\|D_{-}^{\alpha}f\|_{L_2(0, \infty)} \label{7.1} \\
& &  \|f\|_{\widehat{h}_2^{\alpha}}=\|x^{\alpha} (-\frac{1}{x}\frac{d}{dx})^{\alpha}f\|_{L_2(0, \infty)}, \label{7.2}
\end{eqnarray}
here $D_-^{\alpha}$ is the Riemann--Liouville fractional integrodifferentiation, operator in  \eqref{7.2} is defined by
\vskip -10pt
\begin{equation}\label{7.3}{(-\frac{1}{x}\frac{d}{dx})^{\beta}=2^{\beta}I_{-; \, 2, \,0}^{-\beta}x^{-2 \beta},}
\end{equation}
$I_{-; 2, \, 0}^{-\beta}$ is the Erd\'{e}lyi--Kober operator, $\alpha\in\mathbb{R}$.  For $\beta = n \in \mathbb{N}_0$ the
expression \eqref{7.3} reduces to classical derivatives.

\begin{theorem}
Let $f(x) \in \mathbb{D}(0, \infty)$. Then the next formulas are valid:
\begin{eqnarray}
& & D_{-}^{\alpha}f={_1S_{-}^{\alpha-1}} {x^{\alpha} (-\frac{1}{x}\frac{d}{dx})^{\alpha}} f, \label{7.4} \\
& & x^{\alpha} (-\frac{1}{x}\frac{d}{dx})^{\alpha}f={_1P_{-}^{\alpha-1}} D_{-}^{\alpha}f. \label{7.5}
\end{eqnarray}
\end{theorem}

So the Buschman--Erd\'{e}lyi  transmutations of zero order smoothness for $\alpha \in \mathbb{N}$ link differential operators in seminorms definitions \eqref{7.1} and \eqref{7.2}.

\begin{theorem}
Let $f(x) \in \mathbb{D}(0, \infty)$. Then the next inequalities hold true for seminorms:
\vskip -10pt
\begin{eqnarray}
& &  \|f\|_{h_2^{\alpha}} \leq \max (1, \sqrt{1+\sin \pi \alpha}) \|f\|_{\widehat{h}_2^{\alpha}}, \label{7.7}\\
& & \|f\|_{\widehat{h}_2^{\alpha}} \leq \frac{1}{\min (1, \sqrt{1+\sin \pi \alpha})} \|f\|_{h_2^{\alpha}}, \label{7.8}
\end{eqnarray}
here $\alpha$ is any real number except $\alpha \neq -\frac{1}{2}+2k,~k \in \mathbb{Z}$.
\end{theorem}

The constants in inequalities \eqref{7.7}--\eqref{7.8} are not greater than 1, it will be used below.
If $\sin \pi \alpha = -1 $ or  $\alpha = -\frac{1}{2}+2k,~k \in \mathbb{Z}$, then the estimate \eqref{7.8} is not valid.

Define on  $\mathbb{D} (0, \infty )$ the Sobolev norm
\begin{equation}\label{7.9}{\|f\|_{W_2^{\alpha}}=\|f\|_{L_2 (0, \infty)}+\|f\|_{h_2^{\alpha}}.}
\end{equation}
Define one more norm,
\vskip -10pt
\begin{equation}\label{7.10}{\|f\|_{\widehat{W}_2^{\alpha}}=\|f\|_{L_2 (0, \infty)}+\|f\|_{\widehat{h}_2^{\alpha}}}
\end{equation}
Define the spaces $W_2^{\alpha},~ \widehat{W}_2^{\alpha}$ as closures of $D(0,
\infty)$ in  \eqref{7.9} or \eqref{7.10}, respectively.

\begin{theorem}
a) For all $\alpha \in \mathbb{R}$ the space $\widehat{W}_2^{\alpha}$ is continuously imbedded in  $W_2^{\alpha}$, moreover
\begin{equation}\label{7.11}
{\|f\|_{W_2^{\alpha}}\leq A_1 \|f\|_{\widehat{W}_2^{\alpha}},}
\end{equation}
with $A_1=\max (1, \sqrt{1+\sin \pi \alpha})$.

b) Let $\sin \pi \alpha \neq -1$ or $\alpha \neq -\frac{1}{2} + 2k, ~ k \in \mathbb{Z}$.
 Then the inverse embedding of $W_2^{\alpha}$  in $\widehat{W}_2^{\alpha}$ is valid, moreover
\begin{equation}\label{7.12}
{\|f\|_{\widehat{W}_2^{\alpha}}\leq A_2 \|f\|_{W_2^{\alpha}},}
\end{equation}
with $A_2 =1/  \min (1, \sqrt{1+\sin \pi \alpha})$.

c) Let $\sin \pi \alpha \neq -1$, then the spaces $W_2^{\alpha}$  and $\widehat{W}_2^{\alpha}$ are isomorphic with equivalent norms.

d) The constants in embedding inequalities \eqref{7.11}--\eqref{7.12} are  sharp.
\end{theorem}

In fact this theorem is a direct corollary of the results on boundedness and norm estimates in $L_2$ of
the Buschman--Erd\'{e}lyi  transmutations of zero order smoothness. In the same manner, from the unitarity of these operators it follows the next
theorem.

\begin{theorem}
The norms
\begin{eqnarray}
& & \|f\|_{W_2^{\alpha}} = \sum\limits_{j=0}^s \| D_{-}^j f\|_{L_2}, \label{7.13} \\
& & \|f\|_{\widehat{W}_2^{\alpha}}=\sum\limits_{j=0}^s \| x^j(-\frac{1}{x}\frac{d}{dx})^j f \|_{L_2} \label{7.14}
\end{eqnarray}
are equivalent  for integer  $s \in \mathbb{Z}$. Moreover, each term in  \eqref{7.13} equals to appropriate term in \eqref{7.14} of the same index $j$.
\end{theorem}

I.\,Kipriyanov introduced in \cite{Kip2} function spaces which essentially influenced the theory of partial differential equations with Bessel operators and in more general sense, the theory of singular and degenerate equations. These spaces are defined in the following next way.
First we consider subset of even functions in  $\mathbb{D}(0, \infty)$ with all zero derivatives of odd orders at  $x=0$. Denote this set as $\mathbb{D}_c (0, \infty)$ and equip it with a norm
\begin{equation}\label{7.15}{\|f\|_{\widetilde{W}_{2, k}^s} = \|f\|_{L_{2, k}}+\|B_k^{\frac{s}{2}}\|_{L_{2, k}}},
\end{equation}
where $s$ is an even natural number,  $B^{s/2}_k$ is an iteration of the Bessel operator.
Define the Kipriyanov spaces for even  $s$ as a closure of $D_c (0, \infty)$ in the norm \eqref{7.15}.
It is a known fact that equivalent to \eqref{7.15} norm may be defined by \cite{Kip2},
\begin{equation}\label{7.16}
{\|f\|_{\widetilde{W}_{2, k}^s} = \|f\|_{L_{2, k}}+\|x^s(-\frac{1}{x}\frac{d}{dx})^s f\|_{L_{2, k}}}
\end{equation}
So the norm $\widetilde{W}_{2, \, k}^s$ may be defined for all $s$. Essentially this approach is the same as in   \cite{Kip2},
another approach is based on usage of the Hankel transform. Below we adopt the norm \eqref{7.16} for the space $\widetilde{W}_{2, k}^s$.

Define the weighted Sobolev norm by
\begin{equation}\label{7.17}
{\|f\|_{W_{2, k}^s} = \|f\|_{L_{2, k}}+\|D_{-}^s f\|_{L_{2, k}}}
\end{equation}
and a space $W_{2, \, k}^s$ as a closure of $\mathbb{D}_c (0, \infty)$ in this norm.

\begin{theorem}\label{tvloz1}
a) Let $k \neq -n, ~ n \in \mathbb{N}$. Then the space  $\widetilde{W}_{2, \, k}^s$ is continuously embedded into  $W_{2, \, k}^s$, and there exists a constant $A_3>0$ such that
\begin{equation}\label{7.18}{\|f\|_{W_{2, k}^s}\leq A_3 \|f\|_{\widetilde{W}_{2, k}^s},}
\end{equation}
b) Let $k+s \neq -2m_1-1, ~ k-s \neq -2m_2-2, ~ m_1 \in \mathbb{N}_0, ~ m_2 \in \mathbb{N}_0$.
Then the inverse embedding holds true of $W_{2, \, k}^s$ into $\widetilde{W}_{2, \, k}^s$, and there exists a constant $A_4>0$, such that
\begin{equation}\label{7.19}{\|f\|_{\widetilde{W}_{2, k}^s}\leq A_4 \|f\|_{W_{2, k}^s}.}
\end{equation}
c) If the above mentioned conditions are not valid, then the embedding theorems under considerations fail.
\end{theorem}

\begin{corollary}
Let the next conditions hold true: $k \neq -n, ~ n \in \mathbb{N}$; $k+s \neq -2m_1-1,  ~ m_1 \in \mathbb{N}_0; ~ k-s \neq -2m_2-2, ~ m_2 \in \mathbb{N}_0$. Then the Kipriyanov spaces may be defined as closure of $D_c (0, \infty)$ in the weighted Sobolev norm \eqref{7.17}.
\end{corollary}

\begin{corollary}
The sharp constants in embedding theorems \eqref{7.18}--\eqref{7.19} are:
$$
A_3 = \max (1, \|{_1S_-^{s-1}} \| _ {L_{2, k}}), ~ A_4=\max(1, \|{_1P_-^{s-1}}\|_{L_{2, k}}).
$$
\end{corollary}

It is obvious that the theorem above and its corollaries are direct consequences of estimates for
the Buschman--Erd\'{e}lyi  transmutations. The sharp constants in embedding theorems \eqref{7.18}--\eqref{7.19} are also direct  consequences of estimates for the Buschman--Erd\'{e}lyi  transmutations of zero order smoothness.
Estimates in  $L_{p, \alpha}$ not included in this article allow to consider embedding theorems for the general Sobolev and Kipriyanov spaces.

So by applying the Buschman--Erd\'{e}lyi  transmutations of zero order smoothness, we received an answer to a problem which for a long time was discussed in ``folklore": -- are the Kipriyanov spaces isomorphic to power weighted Sobolev spaces or not? Of course we investigated just the simplest case, the results can be generalized to other seminorms, higher dimensions, bounded domains but the principal idea is clear.
All that do not in any sense disparage neither essential role nor necessity for applications of the Kipriyanov spaces in the theory of partial differential equations.

The importance of the Kipriyanov spaces is a special case of the next general principle of L.\,Kudryavtsev:

\begin{center}
 ``\textit{EVERY EQUATION MUST BE INVESTIGATED IN ITS OWN SPACE!}"
\end{center}

The embedding theorems proved in this section  may be applied to direct transfer of known solution estimates for  $B$--elliptic equations in Kipriyanov spaces (cf. \cite{Kip1},\cite{Kip2} ) to new estimates in weighted Sobolev spaces, it is a direct consequence of boundedness and transmutation properties of the Buschman--Erd\'{e}lyi  transmutations.

\smallskip

\subsection{Solution representations to partial differential equations with Bessel operators.} 

The above classes of transmutations may be used for deriving explicit formulas for solutions of partial differential equations with Bessel operators via unperturbed equation solutions. An example is the  $B$--elliptic equation of the form
\begin{equation}
\sum_{k=1}^{n}B_{\nu,x_k}u(x_1,\dots, x_n)=f,
\end{equation}
and similar  $B$--hyperbolic and $B$--parabolic equations. This idea early works by the Sonine--Poisson--Delsarte transmutations, cf. \cite{Car1}--\cite{Car3}, \cite{CarSho}, \cite{Kip1}.
New results follow automatically for new classes of transmutations.

\smallskip

\subsection{Applications to  Dunkl operators.} 

In recent years the Dunkl operators were thoroughly studied. These are difference--differentiation operators consisting of combinations of classical derivatives and finite differences. In higher dimensions, the Dunkl operators are defined by symmetry and reflection groups. For this class there are many results on transmutations which are of Sonine--Poisson--Delsarte and Buschman--Erd\'{e}lyi types, cf. \cite{Sig} and references therein.

\smallskip

\subsection{Applications of Buschman--Erd\'{e}lyi  operators to the Radon transform.} 

It was proved by Ludwig in \cite{Lud} that the Radon transform in terms of spherical harmonics acts in every harmonics at radial components as Buschman--Erd\'{e}lyi  operators. Let us formulate this result.

\begin{theorem}
{\rm (Ludwig theorem, \cite{Lud},\cite{Hel})}
Let the function $f(x)$ be expanded in $\mathbb{R}^n$ by spherical harmonics
\vskip -12pt
\begin{equation}
f(x)=\sum_{k,l}f_{k,l}(r) Y_{k,l}(\theta).
\end{equation}
\vskip -3pt \noindent
Then the Radon transform of this function may be calculated as another series in spherical harmonics,
\vskip -12pt
\begin{equation}
Rf(x)=g(r,\theta)=\sum_{k,l}g_{k,l}(r) Y_{k,l}(\theta),
\end{equation}
\vspace*{-6pt}
\begin{equation}
\label{lu1}
g_{k,l}(r)=ñ(n)\int_r^\infty \lr{1-\frac{s^2}{r^2}}^{\frac{n-3}{2}} C_l^{\frac{n-2}{2}}\lr{\frac{s}{r}}
f_{k,l}(r) r^{n-2}\,ds,
\end{equation}
where $ñ(n)$ is some known constant, $C_l^{\frac{n-2}{2}}\lr{\frac{s}{r}}$ is the Gegenbauer function \cite{BE}.
The inverse formula is also valid of representing values $f_{k,l}(r)$ via $g_{k,l}(r)$.
\end{theorem}

The Gegenbauer function may be easily reduced to the Legendre function,  \cite{BE}.
So the Ludwig formula (\ref{lu1}) reduces the Radon transform in terms of spherical harmonics series and up to unimportant power and constant terms to Buschman--Erd\'{e}lyi  operators of the first kind.

Exactly, this formula in dimension two was developed by Cormack as the first step to the Nobel prize.
Special cases of Ludwig's formula proved in 1966 are for any special spherical harmonics and in the simplest case on pure radial function,
in this case it is reduced to Sonine--Poisson--Delsarte transmutations of Erd\'{e}lyi--Kober type. Besides the fact that such formulas are known  for about half a century they are  rediscovered  still... As consequences of the above connections, the results may be proved for integral representations, norm estimates, inversion formulas for the Radon transform via
Buschman--Erd\'{e}lyi operators. In particular, it makes clear that different kinds of inversion formulas for the Radon transform are at the same time inversion formulas for the Buschman--Erd\'{e}lyi  transmutations of the first kind and vice versa. A useful reference for this approach is \cite{Dea}.

\smallskip

\subsection{Application of the Buschman--Erd\'{e}lyi  transmutations for estimation
of generalized Hardy operators.}

We proved unitarity of the shifted Hardy operators (\ref{3.25}) and mentioned that it is a known fact from \cite{KuPe}.
It is interesting that the Hardy operators naturally arise in transmutation theory. Use Theorem 7 with integer parameter
which guarantees the unitarity for finding more unitary in $L_2(0,\infty)$ integral operators of very simple form.

\smallskip

\begin{theorem}
The next are  pairs of unitary mutually inverse integral operators in $L_2(0,\infty)$:

\newpage 
\begin{eqnarray*}
\label{84}
U_3f= f+\int_0^x f(y)\,\frac{dy}{y},\  U_4f= f+\frac{1}{x}\int_x^\infty f(y)\,dy,\\\nonumber
U_5f= f+3x\int_0^x f(y)\,\frac{dy}{y^2},\  U_6f= f-\frac{3}{x^2}\int_0^x y f(y)\,dy,\\\nonumber
U_7f=f+\frac{3}{x^2}\int_x^\infty y f(y)\,dy,\  U_8f=f-3x \int_x^\infty f(y)\frac{dy}{y^2},\\\nonumber
U_9f=f+\frac{1}{2}\int_0^x \left(\frac{15x^2}{y^3}-\frac{3}{y}\right)f(y)\,dy,\\\nonumber
U_{10}f=f+\frac{1}{2}\int_x^\infty \left(\frac{15y^2}{x^3}-\frac{3}{x}\right)f(y)\,dy.\\\nonumber
\end{eqnarray*}
\end{theorem}


\subsection{Application of the Buschman--Erd\'{e}lyi  transmutations in works of \hfill \break V.\,Katrakhov.}

\medskip

V.\,Katrakhov found a new approach for boundary value problems for elliptic equations with strong singularities of infinite order. For example, for the Poisson equation he studied problems with solutions of arbitrary growth.
At singular point he proposed the new kind of boundary condition: $K$--trace. His results are based on the constant usage of Buschman--Erd\'{e}lyi transmutations of the first kind for definition of norms, solution estimates and correctness proofs \cite{Kat1}--\cite{Kat2}.

\smallskip

Moreover in joint papers with I.\,Kipriyanov he introduced and studied new classes of pseudodifferential operators based on transmutational technics \cite{Kat3}. These results were
paraphrased in reorganized manner also in \cite{Car2}.


\medskip 

\end{document}